\documentclass[11pt, leqno]{article}
\usepackage{amsmath, amsfonts}
\begin{document}
\author{Ajai Choudhry}
\title{ Equal Sums of Like Powers\\ with Minimum Number of Terms}
\date{}
\maketitle

\begin{abstract} 
This paper is concerned with 
 the diophantine system,
$\sum_{i=1}^{s_1} x_i^r=\sum_{i=1}^{s_2} y_i^r,\;
r=1,\,2,\,\ldots,\,k, $
where $s_1$ and $s_2$ are integers such that the total number of terms on both sides, that is, $s_1+s_2,$ is as small as possible. We define $\beta(k)$ to be  the minimum value of $s_1+s_2$ for which there exists a nontrivial solution  of this diophantine system. We find nontrivial integer solutions of this  diophantine system when $k < 6$, and thereby show that $\beta(2) =4,\;\; \beta(3) = 6,\;\; 7 \leq \beta(4) \leq 8$ and $8 \leq \beta(5) \leq 10$. 
\end{abstract}


\section{Introduction}
\hspace*{0.25in}The Tarry-Escott problem of degree $k$  consists of finding two  distinct sets of integers $\{x_1,\,x_2,\,\ldots,\,x_s\}$ and $\{y_1,\,y_2,\,\ldots,\,y_s\}$ such that 
\begin{equation}
\sum_{i=1}^s x_i^r=\sum_{i=1}^s y_i^r,\quad
r=1,\,2,\,\ldots,\,k.
\label{tepsk}
\end{equation}
It is well-known   that for a non-trivial
solution of (\ref{tepsk}) to exist, we must have $s\geq (k+1)$ \cite[p.\  616]{Dor2}. Solutions of (\ref{tepsk}) with the minimum possible value of $s$, that is, with $s=k+1,$ are known as ideal solutions of the problem. 

This paper is concerned with finding solutions in integers of  the related diophantine system,
  \begin{equation}
\sum_{i=1}^{s_1} x_i^r=\sum_{i=1}^{s_2} y_i^r,\quad
r=1,\,2,\,\ldots,\,k, \label{sysgenk} 
\end{equation}
where $s_1$ and $s_2$ are integers such that the total number of terms on both sides, that is, $s_1+s_2,$ is minimum.

 Without loss of generality, we may take  $s_1 \leq  s_2$. A solution of the system of equations \eqref{sysgenk} will be said to be trivial if $y_i=0$ for $s_2-s_1$ values of $i$ and the remaining integers $y_i$ are a permutation of the integers $x_i$. 
 
 We define $\beta(k)$ to be  the minimum value of $s_1+s_2$ for which there exists a nontrivial solution  of the diophantine system \eqref{sysgenk}. 

According to a well-known theorem of  Frolov  \cite[p.\  614]{Dor2},  if $x_i=a_i,\,y_i=b_i, \;i=1,\,2,\,\ldots,\,s$ is any solution of the diophantine system \eqref{tepsk}, and $d$ is an arbitrary integer, then  $x_i=a_i+d,\;y_i=b_i+d,\;i=1,\,2,\,\dots,\,s$, 
is also a solution of \eqref{tepsk}. Taking $d=-a_1$, we immediately get a solution of \eqref{sysgenk} with $s_1=s-1,\;s_2=s$. Thus, if an ideal solution of \eqref{tepsk} is known for  any specific value of $k$,  then $\beta(k) \leq 2k+1.$

Ideal solutions  of \eqref{tepsk} are known for  $k=2,\,3,\,\ldots,\,9$  (\cite{Bor}, \cite{Che}, \cite{Cho1},  \cite{Cho2}, \cite[pp.\  52,\,55-58]{Di1}, \cite{Dor1}, \cite[pp.\  41-54]{Glo}, \cite{Let}, \cite{Smy}) as well as for $k=11$ \cite{Cho3}. Thus, for these values of $k$, we have  $\beta(k) \leq 2k+1$, and in particular, we get
\begin{equation}
\beta(2) \leq 5,\;\;\beta(3) \leq 7,\;\;\beta(4) \leq 9,\;\;\beta(5) \leq 11.
\end{equation}

In this paper, we find parametric solutions of (\ref{sysgenk})  when $2 \leq k \leq 5$, and show that, 
\begin{equation}
\beta(2) =4,\quad \beta(3) = 6,\quad 7 \leq \beta(4) \leq 8,\quad  8 \leq \beta(5) \leq 10. \label{estimatebeta}
\end{equation}

\section{A preliminary lemma}

\noindent {\bf Lemma 1:} If there exists a nontrivial solution of the diophantine system \eqref{sysgenk}, then 
\begin{equation}
{\rm max}(s_1,\,s_2) \geq k+1. \label{betakest1}
\end{equation}
Further, if $k \geq 4$, then,
\begin{equation}
{\rm min}(s_1,\,s_2) \geq 2. \label{betakest2}
\end{equation}

\noindent {\bf Proof:} Any solution of the diophantine system \eqref{sysgenk} with  ${\rm max}(s_1,\,s_2) < k+1$,  implies the existence of a  solution of \eqref{sysgenk} with $s_1=s_2 < k+1$, and by a theorem of Bastien (as quoted by Dickson \cite[p. 712]{Di2}), such a solution is necessarily trivial. Thus we must have the relation \eqref{betakest1}. 

When $k \geq 4$, if we assume a solution of the diophantine system \eqref{sysgenk} with $s_1=1$, it is easily seen  on   eliminating $x_1$ from the two relations obtained by taking $r=2$ and $r=4$ in \eqref{sysgenk} that  the resulting condition  can  be satisfied only if all but one of the numbers $y_i$ are 0, and the solution is necessarily trivial. Thus we must have the relation \eqref{betakest2}.

\noindent {\bf Corollary 1:} For any arbitrary positive integer $k$, 
\begin{equation}
\beta(k) \geq k+2.
 \label{betaestgen}
\end{equation}
Further, when $k \geq 4$, 
\begin{equation}
\beta(k) \geq k+3.
 \label{betaestspl}
\end{equation}

\noindent{\bf Proof:} We have trivially ${\rm min}(s_1,\,s_2)  \geq 1$. The corollary now follows immediately from the lemma.

We have already seen that if there exists a solution of the diophantine system \eqref{sysgenk} with $s_1=s_2$, Frolov's theorem immediately gives another solution of  \eqref{sysgenk} with $s_1=s-1,\;s_2=s$. Thus for  $s_1+s_2$  to be minimum, we  can always take $s_1 < s_2$. Accordingly, we will henceforth always consider the diophantine system \eqref{sysgenk} with $s_1 < s_2$. 

We note that all the equations of the diophantine system \eqref{sysgenk} are homogeneous, and therefore, any  solution of \eqref{sysgenk} in rational numbers may be multiplied through by a suitable constant to obtain a solution of \eqref{sysgenk} in integers.

\section{Determination of $\beta(k)$}
\hspace*{0.25in}It is trivially true that $\beta(1)=3$. In the next four subsections, we will find solutions of the diophantine system \eqref{sysgenk} when $k=2,\,3,\,4$ and $5$ respectively, and prove the results stated in \eqref{estimatebeta}.

\noindent {\bf 3.1}  It follows from Cor. 1 that $\beta(2) \geq 4$. We will show that $\beta(2) = 4$ by solving the diophantine system \eqref{sysgenk} with $s_1=1,\;s_2=3$ and $k=2$. On eliminating 
$x_1$ from the two equations of this diophantine system,   we get,
\begin{equation}
y_1y_2+y_2y_3+y_3y_1=0. \label{sysgen1cond1}
\end{equation}
The complete solution of Eq.~\eqref{sysgen1cond1} is readily obtained and this immediately yields the following simultaneous identities:
\begin{equation}
(p^2+pq+q^2)^r=(p^2+pq)^r+(pq+q^2)^r+(-pq)^r,\;\;r=1,\,2,
 \end{equation}
 where $p$ and $q$ are arbitrary parameters. This shows that $\beta(2) = 4$.

\noindent {\bf 3.2} We now consider the diophantine system \eqref{sysgenk} when $k=3$. It follows from Cor. 1 that $\beta(3) \geq 5$. We will prove that $\beta(3)=6$.

If there exists a  nontrivial solution of the system of equations,
\begin{equation}
\sum_{i=1}^{s_1} x_i^r=\sum_{i=1}^{s_2} y_i^r,\quad
r=1,\,2,\,3, \label{sysgen3} 
\end{equation}
with $s_1 < s_2$, it follows from Lemma 1 that $s_2 \geq 4$. Thus   the only case to consider when $s_1+s_2=5$ is with $s_1=1$ and $s_2=4$ when  the diophantine system \eqref{sysgen3}  may be written as,
\begin{align}
x_1&=y_1+y_2+y_3+y_4, \label{tep3r1}\\
x_1^2&=y_1^2+y_2^2+y_3^2+y_4^2, \label{tep3r2} \\
x_1^3&=y_1^3+y_2^3+y_3^3+y_4^3. \label{tep3r3}
\end{align}

Eliminating $x_1$ from  Eqs.~\eqref{tep3r1} and \eqref{tep3r2}, we get the condition,
\begin{equation}
(y_1+y_2+y_3)y_4+y_1y_2+y_1y_3+y_2y_3=0,\label{restep3r1r2}
\end{equation}
while eliminating $x_1$ from Eqs.~\eqref{tep3r1} and \eqref{tep3r3}, we get the condition,
\begin{equation}
(y_1+y_2+y_3)y_4^2+(y_1+y_2+y_3)^2y_4+(y_2+y_3)(y_1+y_3)(y_1+y_2)=0. \label{restep3r1r3}
\end{equation}

If $y_1+y_2+y_3=0$, it follows from \eqref{tep3r1} that $x_1=y_4$ and now \eqref{tep3r2} implies that $y_1=0,\,y_2=0,\,y_3=0$, and we thus get a trivial solution of the system of equations  \eqref{tep3r1}, \eqref{tep3r2} and \eqref{tep3r3}. If $y_1+y_2+y_3 \neq 0$, we eliminate $y_4$ from Eqs.~\eqref{restep3r1r2} and \eqref{restep3r1r3}, and get,
\begin{equation}
(y_1^2+y_1y_2+y_2^2)y_3^2+y_1y_2(y_1+y_2)y_3+y_1^2y_2^2=0. \label{restepr123}
\end{equation}

Eq.~\eqref{restepr123} can have a rational solution for $y_3$ only if its discriminant, that is, $-(3y_1^2+2y_1y_2+3y_2^2)y_1^2y_2^2$, is  a perfect square. It follows that either $y_1$ or $y_2$, or both of them, must be 0, and in each case, it readily follows that the only solution of  Eqs.~\eqref{tep3r1}, \eqref{tep3r2}  and \eqref{tep3r3} is the trivial solution. 

Since the diophantine system \eqref{sysgen3} has only the trivial solution when  $s_1+s_2 = 5$,  it follows that, 
\begin{equation}  
\beta(3) > 5. \label{beta3est1}
\end{equation}

We will now  obtain nontrivial solutions of the diophantine system \eqref{sysgen3} with $s_1=2$ and $s_2=4$, that is, of the system of equations,
\begin{align}
x_1+x_2&=y_1+y_2+y_3+y_4,\label{sysgen3eq1} \\
x_1^2+x_2^2&=y_1^2+y_2^2+y_3^2+y_4^2,\label{sysgen3eq2}\\
x_1^3+x_2^3&=y_1^3+y_2^3+y_3^3+y_4^3.\label{sysgen3eq3}
\end{align}

If $(a,\,b,\,c)$ is any Pythagorean triple satisfying the relation $a^2+b^2=c^2$, it is easily seen that a  solution in integers of the  simultaneous equations  \eqref{sysgen3eq1}, \eqref{sysgen3eq2} and \eqref{sysgen3eq3} is given by  $(x_1,\,x_2,\,y_1,\,y_2,\,y_3,\,y_4)=(c,\,-c,\,a,\,-a,\,b,\,-b)$. 

Next we obtain a more general parametric  solution of  the simultaneous equations \eqref{sysgen3eq1}, \eqref{sysgen3eq2} and \eqref{sysgen3eq3}. We will use a  parametric solution of the simultaneous diophantine  equations,
\begin{equation}
\begin{aligned}
x+y+z&=u+v+w,\\
x^3+y^3+z^3&=u^3+v^3+w^3,
\end{aligned}
\end{equation}
given in \cite[Theorem 1, p.\ 61]{Cho3}. From this solution, on writing $x_1=x,\,x_2=y,\,y_1=u,\,y_2=v,\,y_3=w,\,y_4=-z$, we immediately derive the following solution of the simultaneous equations \eqref{sysgen3eq1} and \eqref{sysgen3eq3} in terms of arbitrary parameters $p,\,q,\,r$ and $s$:
\begin{equation}
\begin{aligned}
x_1&= pq-pr+qr-(p-q-r)s,\\
x_2& = -pq+pr+qr+(p-q+r)s, \\
y_1& = pq+pr-qr+(p-q+r)s, \\
 y_2& = pq-pr+qr+(p+q-r)s, \\
y_3& = -pq+pr+qr-(p-q-r)s,\\
y_4& = -pq-pr+qr-(p+q-r)s.
\end{aligned}
\label{valpqrssysgen3}
\end{equation}

Substituting the above values of $x_i,\;y_i$ in \eqref{sysgen3eq2}, we get, after necessary transpositions, the following quadratic equation in $s$:
\begin{multline}
 2(p+q-r)^2s^2+(12p^2q-4p^2r-4pq^2-4pqr\\
 +4pr^2+4q^2r-4qr^2)s+2(pq+pr-qr)^2=0. \label{sysgen3eq2a}
\end{multline}
On taking $r=p+q$, the coefficient of $s^2$ in Eq,~\eqref{sysgen3eq2a} vanishes, and we can then readily  solve Eq.~\eqref{sysgen3eq2a} for $s$, and we thus obtain the following solution of the simultaneous equations \eqref{sysgen3eq1}, \eqref{sysgen3eq2} and \eqref{sysgen3eq3} in terms of arbitrary parameters $p$ and $q$:
\begin{equation}
\begin{aligned}
x_1 &= (3p^4-2p^3q-p^2q^2+q^4)q,\;\; &x_2 &= (p^4-p^2q^2-2pq^3+3q^4)p, \\
y_1 &= (p^4-p^2q^2+2pq^3-q^4)p, &y_2 &= 2pq(p-q)(p^2-pq-q^2), \\
y_3 &= -(p^4-2p^3q+p^2q^2-q^4)q, &y_4 &= 2pq(p-q)(p^2+pq-q^2).
\end{aligned}
\label{solsysgen3}
\end{equation}

As a numerical example, taking $p=2,\;q=1$, we get the solution,
\[29^r+22^r=30^r+4^r+(-3)^r+20^r,\quad r=1,\,2,\,3.
\]

We note that more parametric solutions of the system of equations \eqref{sysgen3eq1}, \eqref{sysgen3eq2} and \eqref{sysgen3eq3} can be obtained by solving Eq.~\eqref{sysgen3eq2a}  in different ways, for example, by choosing $p, \,q,\,r$ such that the constant term in Eq.~\eqref{sysgen3eq2a} vanishes and then solving this equation for $s$, or by  choosing $p, \,q,\,r$ such that the discriminant of  Eq.~\eqref{sysgen3eq2a}, considered as a quadratic equation in $s$, becomes a perfect square,   and then solving this equation for $s$. 

As we have obtained nontrivial solutions of the system of equations \eqref{sysgen3eq1}, \eqref{sysgen3eq2} and \eqref{sysgen3eq3}, we get $\beta(3) \leq 6$, and on combining this result with \eqref{beta3est1}, we get,
\begin{equation}
\beta(3)=6. \label{valbeta3}
\end{equation}

\noindent {\bf 3.3} We will  now obtain parametric solutions of the diophantine system,
\begin{equation}
\sum_{i=1}^3 x_i^r=\sum_{i=1}^5 y_i^r,\quad
r=1,\,2,\,3,\,4. \label{sysgen4} 
\end{equation}

We write,
\begin{equation}
\begin{aligned}
x_1 &= 4uv+w+1,& x_2 &= -4uv+w-1,\\
 x_3 &= -8u^2+8uv+4u-2,& y_1 &= 4u-2,\\
  y_2 &= -4u,& y_3 &= 4uv+w-1,\\
   y_4 &= -4uv+w+1,& y_5 &= -8u^2+8uv+4u,
   \end{aligned}
   \label{sysgen4valxy}
\end{equation}
where $u,\,v,\,w$ are arbitrary parameters. 

It is readily verified that  the values of $x_i,\,y_i$ given by \eqref{sysgen4valxy} satisfy Eq.~\eqref{sysgen4} when $r=1$ and $r=2$. Further, on substituting these values of $x_i,\,y_i$ in Eq.~\eqref{sysgen4} and taking $r=3$, we get the condition,
\begin{equation}
4u^3-8u^2v+4uv^2-4u^2+4uv-vw+u-v=0. \label{sysgen4condr3}
\end{equation}
On solving Eq.~\eqref{sysgen4condr3}, we get,
\begin{equation}
w=(4u^3-8u^2v+4uv^2-4u^2+4uv+u-v)/v.\label{sysgen4valw}
\end{equation}

Finally, we substitute the values of $x_i,\,y_i$ given by \eqref{sysgen4valxy} in Eq.~\eqref{sysgen4}, and take $r=4$, and use  the value of $w$ given by \eqref{sysgen4valw} to get the condition,
\begin{equation}
u^2(2u-1)^2\{24uv^2-2(4u-1)^2v+3u(2u-1)^2\}=0. \label{sysgen4condr4}
\end{equation}

While equating the first two factors of Eq.~\eqref{sysgen4condr4} to 0 leads to trivial results, on equating the last factor to 0, we get a quadratic equation in $v$ which will have a rational solution if its discriminant $4(-32u^4+32u^3+24u^2-16u+1)$ becomes a perfect square. We thus have to solve the diophantine equation,
\begin{equation}
t^2=-32u^4+32u^3+24u^2-16u+1. \label{ecquartic}
\end{equation}

Now Eq.~\eqref{ecquartic} is a quartic model of an elliptic curve, and we use  the birational transformation given by,
\begin{equation}
\begin{aligned}
t &= (X^3-36X^2+36X-72Y+432)/(4X+Y-12)^2,\\
 u& = (X-12)/(4X+Y-12),
 \end{aligned}
\label{birat1}
\end{equation}
and,
\begin{equation}
\begin{aligned}
X &=(4u^2-8u+t+1)/(2u^2),\\
 Y& = (8u^3+12u^2-4ut-12u+t+1)/(2u^3),
 \end{aligned}
 \end{equation}
 to reduce Eq.~\eqref{ecquartic} to the Weierstrass model which is given by the cubic equation,
 \begin{equation}
 Y^2=X^3-36X. \label{ecweier}
 \end{equation}

It is readily seen from Cremona's well-known tables \cite{Cre} 
that \eqref{ecweier} is an elliptic curve of rank 1 and its Mordell-Weil basis is given  by the rational point $P$ with co-ordinates $(X,\,Y)=(-3,\,9).$ There are thus infinitely many rational points on the elliptic curve \eqref{ecweier} and these can be obtained by the group law. Using the relations \eqref{birat1}, we can find infinitely many rational solutions of Eq.~\eqref{ecquartic} and thus obtain infinitely many integer solutions of the diophantine system \eqref{sysgen4}.

While the point $P$ leads to a trivial solution of the diophantine system \eqref{sysgen4}, the point $2P$ yields the solution,
\[
(-74)^r+124^r+78^r=126^r+(-72)^r+(-20)^r+70^r+24^r,\;\;r=1,\,2,\,3,\,4,
\]
and the point $3P$ leads to the solution,
\begin{multline*}
(-40573)^r+66494^r+118981^r=(-15181)^r+119510^r+63756^r\\
+(-37835)^r+14652^r,\;\;r=1,\,2,\,3,\,4.
\end{multline*}

In view of the above solutions of the diophantine system \eqref{sysgen4}, it follows that $\beta(4) \leq 8$, and on combining with the result $\beta(4) \geq 7$ which follows from Cor. 1, we get,
\[
7 \leq \beta(4) \leq 8. 
\]


\noindent {\bf 3.4} We will  now obtain parametric solutions of the diophantine system,
\begin{equation}
\sum_{i=1}^4 x_i^r=\sum_{i=1}^6 y_i^r,\quad
r=1,\,2,\,3,\,4,\,5. \label{sysgen5} 
\end{equation}

We will  use   a parametric solution of the diophantine system,
\begin{equation}
\sum_{i=1}^{6} X_i^r=\sum_{i=1}^{6} Y_i^r,\quad
r=1,\,2,\,3,\,4,\,5, \label{tep5} 
\end{equation}
 to obtain two parametric solutions of the diophantine system \eqref{sysgen5}. 
 
A solution of  the simultaneous equations $\sum_{i=1}^{3} X_i^r=\sum_{i=1}^{3} Y_i^r,\;r=2,\,4$, in terms of arbitrary parameters $m,\,n,\,x$ and $y$, given in 
\cite[p. 102]{Cho2}, is as follows:
\begin{equation}
\begin{aligned}
X_1&=(m+2n)x-(m-n)y,\quad &X_2&=-(2m+n)x-(m+2n)y,\\
X_3&=(m-n)x+(2m+n)y,&Y_1&=(m-n)x-(m+2n)y,\\
Y_2&=-(2m+n)x-(m-n)y,&Y_3&=(m+2n)x+(2m+n)y,
\end{aligned}
\label{soltep5}
\end{equation}
It immediately follows that a parametric solution of the diophantine system \eqref{tep5} is given by \eqref{soltep5} and
\begin{equation}
\begin{aligned}
X_4&=-X_3,\;\;&X_5&=-X_2,\;\;&X_6&=-X_1,\\
Y_4&=-Y_3,&Y_5&=-Y_2,&Y_6&=-Y_1.
\end{aligned}
\label{valX456}
\end{equation}

We choose the parameters $x$ and $y$ such that $X_3=0$ and immediately obtain the following parametric solution of the diophantine system \eqref{sysgen5}: 
\begin{equation}
\begin{aligned}
x_1 &= m^2+mn+n^2,\quad & x_2 &= m^2+mn+n^2,\\ x_3 &= -m^2-mn-n^2, &x_4 &= -m^2-mn-n^2,\\ y_1 &= m^2-n^2, &y_2 &= -m^2-2mn,\\ y_3 &= 2mn+n^2,& y_4 &= -2mn-n^2,\\ y_5 &= m^2+2mn, &y_6 &= -m^2+n^2, 
\end{aligned}
\label{sol1sysgen5}
\end{equation}
where $m$ and $n$ are arbitrary parameters.

To obtain  a second solution of the diophantine system \eqref{sysgen5}, we again use the  parametric solution of the diophantine system \eqref{tep5} given by \eqref{soltep5} and \eqref{valX456}. We now  choose the parameters $x,\,y$ such that we get $X_2=X_3$, and then  apply the theorem of Frolov mentioned in the Introduction, taking $d=-X_3$. We thus get a solution of the diophantine system \eqref{sysgen5} which may be written as follows: 
\begin{equation}
\begin{aligned}
x_1 &= 3m^2+3mn+3n^2,& x_2 &= 2m^2+2mn+2n^2,\\ x_3 &= -m^2-mn-n^2,& x_4 &= 2m^2+2mn+2n^2,\\ y_1 &= 3m^2+3mn,& y_2 &= -3mn,\\ y_3 &= 3mn+3n^2,& y_4 &= 2m^2-mn-n^2,\\ y_5 &= 2m^2+5mn+2n^2,& y_6 &= -m^2-mn+2n^2,
\end{aligned}
\label{sol2sysgen5}
\end{equation}
where $m$ and $n$ are arbitrary parameters.

As a numerical example, taking $m=2,n=1$ in \eqref{sol2sysgen5}, we get the solution,
\[
21^r+14^r+(-7)^r+14^r=18^r+(-6)^r+9^r+5^r+20^r+(-4)^r,\;\;r=1,\,2,\,3,\,4,\,5.
\]

The two  parametric solutions \eqref{sol1sysgen5} and \eqref{sol2sysgen5} of the diophantine system \eqref{sysgen5} are rather special since in both of them, the ratios $x_i/x_j$ of the integers on the left-hand side are all fixed. We now show that there exist infinitely many other  solutions of the diophantine system \eqref{sysgen5} that are not generated by these parametric solutions. 

We write,
\begin{equation}
\begin{aligned}
x_1&=uv^2+(6u^3-12u^2+32u-32)v\\
& \quad \quad +9u^5-36u^4-336u^2+96u^3+240u,\\
 x_2& = (2u-2)v^2+(12u^3-48u^2+40u-16)v\\
& \quad \quad +18u^5-126u^4-288u^2+264u^3+96,\\
x_3&=-x_1,\quad x_4=-x_2,\\
  y_1 &= (2u-2)v^2+(12u^3-48u^2+48u)v\\
& \quad \quad +18u^5-126u^4-144u^2+288u^3+96u-96,\\
   y_2 &= uv^2+(6u^3-12u^2-32u+32)v\\
& \quad \quad +9u^5-36u^4+240u^2-96u^3-144u,\\
    y_3 &= 2v^2+(24u^2-40u+16)v\\
& \quad \quad +54u^4-264u^3+96u+192u^2-96,\\
y_4&=-y_3,\quad y_5=-y_2,\,\quad y_6=-y_3.
    \end{aligned}
    \label{valxysysgen5}
    \end{equation}

 With these values of $x_i,\,y_i$, it is readily seen  that \eqref{sysgen5} is identically satisfied for $r=1,\,3,\,5$. Further,
 \begin{align}
  \sum_{i=1}^4 x_i^2-\sum_{i=1}^6 y_i^2&=-8(9u^4-72u^3+24u^2+96u-48-v^2)^2,\\
  \sum_{i=1}^4 x_i^4-\sum_{i=1}^6 y_i^4&=-8(9u^4-72u^3+24u^2+96u-48-v^2)^4.
  \end{align}

It follows that a solution of the diophantine system \eqref{sysgen5} will be given by \eqref{valxysysgen5} if we choose $u,\,v$ such that,
\begin{equation}
v^2=9u^4-72u^3+24u^2+96u-48. \label{sysgen5ecquartic}
\end{equation}
Now Eq.~\eqref{sysgen5ecquartic} represents the  quartic model of an elliptic curve, and the birational transformation given by,
\begin{equation}
\begin{aligned}
u& = (6X+2Y-12)/(3X-24),\\
 v& = (4X^3-96X^2+84X-144Y+832)/\{3(X-8)^2\},
 \end{aligned} \label{sysgen5birat1}
 \end{equation}
and,
\begin{equation}
\begin{aligned}
X&=(9u^2-36u+3v+4)/8,\\
Y&=(27u^3-162u^2+9uv+36u-18v+72)/16,
\end{aligned}\label{sysgen5birat2}
\end{equation}
reduces Eq.~\eqref{sysgen5ecquartic} to the Weierstrass form of the elliptic curve which is as follows:
\begin{equation}
Y^2=X^3-21X-20. \label{sysgen5ecweier}
\end{equation}

We again refer to Cremona's  database of elliptic curves,  
and find that  \eqref{sysgen5ecweier} represents an elliptic curve of rank 1 and its Mordell-Weil basis is given  by the rational point $P$ with co-ordinates $(X,\,Y)=(-3,\,4).$ There are thus infinitely many rational points on the elliptic curve \eqref{ecweier} and these can be obtained by the group law. Using the relations \eqref{sysgen5birat1}, we can find infinitely many rational solutions of Eq.~\eqref{sysgen5ecquartic} and thus obtain infinitely many solutions of the diophantine system \eqref{sysgen5}. 
While the point P leads to a trivial solution of the diophantine system \eqref{sysgen5}, the point 2P yields the solution, 
\begin{multline*}
241^r+218^r+(-241)^r+(-218)^r=266^r+143^r+120^r\\
+(-266)^r+(-143)^r+(-120)^r,\;\;r=1,\,2,\,3,\,4,\,5.
\end{multline*}

 
 The solutions of the diophantine system \eqref{sysgen5} obtained in this Section show that $\beta(5) \leq 10$, and on combining with the result $\beta(5) \geq 8$ which follows from Cor. 1, we get,
 \[ 8 \leq \beta(5) \leq 10.
 \]
 
 \section{Concluding Remarks} 
 \hspace*{0.25in} It would be of interest to determine the precise values of $\beta(4)$ and $\beta(5)$. Further, it would be interesting to find integer solutions of the diophantine system \eqref{sysgenk} with  $k>5$ and $s_1+s_2 < 2k+1$.

\begin{center}
\Large
Acknowledgment
\end{center}
 
I wish to  thank the Harish-Chandra Research Institute, Allahabad for providing me with all necessary facilities that have helped me to pursue my research work in mathematics.

\medskip

\noindent Postal Address: Ajai Choudhry, \\
\noindent \hspace*{1.05 in} 13/4 A Clay Square, \\
\noindent \hspace*{1.05 in} Lucknow-226001, INDIA.

\noindent E-mail address: ajaic203@yahoo.com

\end{document}